\documentclass[12pt]{amsart}

\usepackage{amsthm}
\usepackage{amsmath}
\usepackage{amssymb}
\usepackage{eucal}
\usepackage{amsfonts}
\usepackage{amsmath}
\usepackage[Lenny]{fncychap}
\usepackage{enumerate}
\usepackage{amssymb}
\usepackage[english]{babel}
\usepackage[breaklinks=true]{hyperref}
\usepackage[all]{xy}
\usepackage{fancybox}
\usepackage{latexsym,mathrsfs,longtable,lscape,dsfont,yfonts}
\usepackage{srcltx}
\usepackage{soul}

\usepackage[normalem]{ulem}

\ifx\pdfpagewidth\undefined 
\usepackage[T1]{fontenc}
\else 
\usepackage[OT1]{fontenc} 
\fi 
\usepackage{amsfonts,amssymb,latexsym,longtable,amsmath}

\newtheorem{thm}{Theorem}[section]
\newcommand{\bthm}{\begin{thm}} \newcommand{\ethm}{\end{thm}}
\newtheorem{prop}[thm]{Proposition}
\newcommand{\bprp}{\begin{prop}} \newcommand{\eprp}{\end{prop}}
\newtheorem{fact}[thm]{Fact}
\newcommand{\bfct}{\begin{fact}} \newcommand{\efct}{\end{fact}}
\newtheorem{prob}[thm]{Problem}
\newcommand{\bprb}{\begin{prob}} \newcommand{\eprb}{\end{prob}}
\newtheorem{quest}[thm]{Question}
\newcommand{\bqtn}{\begin{quest}} \newcommand{\eqtn}{\end{quest}}
\newtheorem{lem}[thm]{Lemma}
\newcommand{\blem}{\begin{lem}} \newcommand{\elem}{\end{lem}}
\newtheorem{claim}[thm]{Claim}
\newcommand{\bclm}{\begin{claim}} \newcommand{\eclm}{\end{claim}}
\newtheorem{cor}[thm]{Corollary}
\newcommand{\bcor}{\begin{cor}} \newcommand{\ecor}{\end{cor}}
\newtheorem{conj}[thm]{Conjecture}
\newcommand{\bcnj}{\begin{conj}} \newcommand{\ecnj}{\end{conj}}
\theoremstyle{definition}
\newtheorem{defn}[thm]{Definition}
\newcommand{\bdfn}{\begin{defn}} \newcommand{\edfn}{\end{defn}}
\newtheorem{spec}[thm]{Specializing}
\newcommand{\bspc}{\begin{spec}} \newcommand{\espc}{\end{spec}}
\theoremstyle{remark}
\newtheorem{rem}[thm]{Remark}
\newcommand{\brem}{\begin{rem}} \newcommand{\erem}{\end{rem}}
\newtheorem{cnv}[thm]{Convention}
\newcommand{\bcnv}{\begin{cnv}} \newcommand{\ecnv}{\end{cnv}}
\newtheorem{exam}[thm]{Example}
\newcommand{\bexm}{\begin{exam}} \newcommand{\eexm}{\end{exam}}
\newcommand{\bpf}{\begin{proof}} \newcommand{\epf}{\end{proof}}

\newtheorem{thmy}{\textbf{Theorem}}
\newenvironment{thmx}{\stepcounter{thm}\begin{thmy}}{\end{thmy}}

\newcommand{\C}{\mathbb C}

\newcommand{\T}{\mathbb T}
\newcommand{\N}{\mathbb N}

\renewcommand{\phi}{\varphi}
\renewcommand{\theta}{\vartheta}

\newcommand{\gep}{{\epsilon}}

\newcommand{\U}{\mathbb{U}}

\newcommand{\mkp}{\medskip}
\newcommand{\bkp}{\bigskip}
\def\defi{\buildrel\rm def \over=}

\begin{document}

\title[Interpolation sets]{Interpolation sets in spaces of continuous metric-valued functions}

\author[M. Ferrer]{Mar\'ia V. Ferrer}
\address{Universitat Jaume I, Instituto de Matem\'aticas de Castell\'on,
Campus de Riu Sec, 12071 Castell\'{o}n, Spain.}
\email{mferrer@mat.uji.es}

\author[S. Hern\'andez]{Salvador Hern\'andez}
\address{Universitat Jaume I, Departamento de Matem\'{a}ticas,
Campus de Riu Sec, 12071 Castell\'{o}n, Spain.}
\email{hernande@mat.uji.es}

\author[L. T\'arrega]{Luis T\'arrega}
\address{Universitat Jaume I, IMAC and Departamento de Matem\'{a}ticas,
Campus de Riu Sec, 12071 Castell\'{o}n, Spain.}
\email{ltarrega@uji.es}

\thanks{ Research Partially supported by the Spanish Ministerio de Econom\'{i}a y Competitividad, grant MTM2016-77143-P,
and the Universitat Jaume I, grant P1171B2015-77 (AEI/FEDER, UE). The second
author also acknowledges partial support by Generalitat Valenciana,
grant code: PROMETEO/2014/062.}

\begin{abstract}
Let $X$ and $K$ be a \v{C}ech-complete topological group and a compact group, respectively. We prove that
if $G$ is a non-equicontinuous subset of $CHom(X,K)$, the set of all continuous homomorphisms of
$X$ into $K$, then there is a countably infinite subset $L\subseteq G$ such that $\overline L^{K^X}$ is canonically
homeomorphic to $\beta\omega$, the Stone-\v{C}ech compactifcation of the natural numbers. As a consequence,
if $G$ is an infinite subset of $CHom(X,K)$ such that for every countable subset $L\subseteq G$
and compact separable subset $Y\subseteq X$ it holds that either
$\overline{L}^{K^Y}$ has countable tightness or $\vert\overline{L}^{K^Y}\vert \leq \mathfrak{c} $, then $G$ is equicontinuous.
Given a topological group $G$, denote by $G^+$ the (algebraic) group $G$ equipped with the Bohr topology.
It is said that $G$ \emph{respects} a topological property $\mathcal P$ when $G$ and $G^+$ have the same subsets satisfying $\mathcal P$.
As an application of our main result, we prove that if $G$ is an abelian, locally quasiconvex, locally $k_\omega$ group, then the following holds:
(i) $G$ respects any \emph{compact-like} property $\mathcal P$ stronger than or equal to functional boundedness; (ii)
$G$ strongly respects compactness.
\end{abstract}

\thanks{{\em 2010 Mathematics Subject Classification.} Primary 22A05; 43A46; 54H11. Secondary 22D35; 43A40; 54C45; 46A99\\
{\em Key Words and Phrases:} \v{C}ech-complete group, Locally $k_{\omega}$-group; Pontryagin duality; $I_0$-set; Interpolation set;
Bohr compactification; Bohr topology; respects compactness}


\date{\today}

\maketitle \setlength{\baselineskip}{24pt}

\section{Introduction}


Let $G$ be a locally compact abelian group and let $\widehat G$ denote its Pontryagin dual group. We say that a subset
$E$ of $\widehat G$ is \emph{Sidon} if for every bounded function $f$ on $E$ there corresponds a Borel regular measure on $G$, $\mu$,
such that $\widehat \mu(\gamma)=f(\gamma)$ for all $\gamma\in E$ (here $\widehat \mu$ denotes the Fourier transform of $\mu$).
If, in addition, $\mu$ is assumed to be discrete (it has a countable support) then it is said that $E$ is an $I_0$-\emph{set}.
Therefore, each $I_0$-set is Sidon. For instance, lacunary (or Hadamard) sets (i.e. sequences $(z_n)_n\subseteq \N$ such that
$\inf z_{n+1}/z_n> 1 $) are perhaps the simplest examples of $I_0$-sets.
The search for interpolation sets is a main goal in harmonic analysis and the monograph by Graham and Hare \cite{Graham2013} contains
most of the recent results in this area.

In this paper, this question is approached from a topological viewpoint that is based on the equivalent formulation of this notion given
by Hartman and Ryll-Nardzewski \cite{Hartman1964} (in fact, the term $I_0$-set is due to them).
According to their (equivalent) definition a subset $E$ of a locally compact abelian group $G$ is an $I_0$-set if for each $f\in l^\infty(E)$
(that is, for each complex-valued, bounded function defined on $E$) there exists an almost periodic function
$f^b$ on $G$ such that $f(\gamma)=f^b(\gamma)$ for all $\gamma\in E$. Furthermore, since every almost periodic function on a topological group $G$
is the restriction of a continuous function defined on the Bohr compactification $bG$ of $G$, it follows that $E\subseteq G$ is
an $I_0$-set if each $f\in l^\infty(E)$ can be extended to a continuous function $f^b$ defined on $bG$. The latter property
implies that $\overline E^{bG}$ is canonically homeomorphic to $\beta E$, the Stone-\v{C}ech compactification of the set $E$
equipped with the discrete topology. This equivalent definition of $I_0$-set
and the \emph{duality methods} obtained from Pontryagin-van Kampen duality allows us to apply topological techniques in the investigation
of this family of interpolation sets. Thus, we can prove the existence of $I_0$-sets for much larger classes of groups  than
locally compact abelian groups. Several applications of our results to different questions related to the Bohr compactification and topology of
topological abelian groups are also obtained. Last but not least, we deal with the topological properties of sets of continuous functions.
Indeed, if $X$ and $M$ are a topological space and a metrizable space respectively, given a subset $G\subseteq C(X,M)$, we look at the possible
existence of copies of $\beta\omega$ (the Stone-\v{C}ech compactification of the natural numbers) within $\overline{G}^{M^X}$.
This property, or its absence, has deep implications on the topological structure of $G$ as a set of continuous functions on $X$
and has found many applications in different settings (for instance, see \cite{Galindo_Hernandez2004,Glasner2012,Filali2013,Graham2013}
where there are applications to topological groups, dynamical systems, functional analysis and harmonic analysis, respectively).


The starting point of this paper stems from a celebrated theorem by Bourgain, Fremlin and Talagrand
about compact subsets of Baire class $1$ functions \cite{Bourgain1978a}, that we present in the way
it is formulated by Todor\v{c}evi\'{c} in \cite{Todorcevic1997}.

\bthm\emph{(J. Bourgain, D.H. Fremlin, M. Talagrand)}\label{bft}
Let $X$ be a Polish space and let $\lbrace f_n\rbrace_{n<\omega}\subseteq C(X)$ be a pointwise bounded sequence.
The following assertions are equivalent (where the closure is taken in $\mathbb R^X$):
\begin{enumerate}[(a)]
\item $\lbrace f_n\rbrace_{n<\omega}$ is sequentially dense in its closure.
\item The closure of $\lbrace f_n\rbrace_{n<\omega}$ contains no copy of $\beta\omega$.
\end{enumerate}
\ethm
\mkp

A variant of this result is due to Pol \cite[p. 34]{Pol1984}, that again was formulated in different terms (cf. \cite{Cascales2000}).
Here $B_1(X)$ denotes the set of all Baire class 1 functions defined on $X$.

\bthm\emph{(R. Pol)}\label{pol}
Let $X$ be a complete metric space, $G$ a subset of $C(X)$, which is uniformly bounded, and $K = \overline{G}$ the closure of $G$ in $[-1,1]^X$.
Then the following are equivalent:
\begin{enumerate}[(a)]
\item $\overline G^{\mathbb R^X}\nsubseteq B_1(X)$.
\item $G$ contains a sequence whose closure in $\mathbb R^X$ is homeomorphic to $\beta\omega$.
\end{enumerate}
\ethm
\mkp

In both cases we have a dichotomy result that basically characterizes two crucial properties about sets of continuous functions defined on a Polish
and metric complete space, respectively.
In this paper we look at this question 
in terms of the set of continuous functions $G\subseteq C(X,M)$ alone, when $M$ is a general metric space.
We first extend the notion of $I_0$-set, given by Hartman and Ryll-Nardzewski for complex-valued functions,
to a more general setting, which will be needed later on when we apply it to topological groups.

\bdfn \label{def_InterpolationSet}
Let $X$ and $M$ be a topological space and metric space, respectively. If $C(X,M)$ denotes the set of all continuous functions from $X$ to $M$,
we say that a subset $Y$ of $X$ is an \emph{$M$-interpolation set} (or, we can simply say an \emph{Interpolation set} for $C(X,M)$)
when for each function $g\in M^Y$, which has relatively compact range in $M$, there exists a map $f\in C(X,M)$ such that $f_{|Y}=g$.
\edfn \mkp

A main goal in this paper is the understanding of the key (topological) facts that characterize the existence of interpolation sets.
Thereby, this research continues the task accomplished in previous projects \cite{Galindo_Hernandez1999,Galindo_Hernandez2004} and \cite{Fer_Her_Tar_dichotomy}.
Here, 
we introduce a crucial property stronger than the mere \emph{non-equicontinuity}, that provides sufficient conditions for the existence of
Interpolation sets in different settings. We refer to \cite{Bourgain1977} for its motivation, where this notion implicitly appears.

\bdfn\label{def_Bourgain}
Let $X$ be a topological space and let $M$ be a metric space. We say that $G\subseteq C(X,M)$
is a $\frak{B}$-family if the following two conditions hold:
\begin{enumerate}[(a)]
\item $G$ is relatively compact in $M^X$.
\item There exists a nonempty open set $V$ of $X$ and $\epsilon>0$ such that for every finite collection
$\lbrace U_1,\ldots,U_n\rbrace$ of nonempty relatively open sets of $V$ there is a $g\in G$ such that
$\hbox{diam}(g(U_{j}))\geq\epsilon$ for all $j\in \lbrace 1,\ldots,n\rbrace$.\\
\end{enumerate}
\edfn

\brem
In \cite{Fer_Her_Tar2017}, we defined a subset $G$ of $C(X,M)$ as \emph{almost equicontinuous}
(resp. \emph{hereditarily almost equicontinuous})
if $G$ is equicontinuous on a dense subset of $X$ (resp. if $G$ is almost equicontinuous
for every closed nonempty subset of $X$). We do not know which is the relation between
the notions of being a $\frak{B}$-family and the negation of being almost equicontinuous
or hereditary almost equicontinuous when $X$ is a \v Cech-complete space.
However, in the cases in which this relation is known (topological groups, for instance),
the existence of interpolation sets is assured as we show below.
\erem \mkp


\bdfn
A map $f\colon X\to Y$ defined between two topological spaces $X$ and $Y$ is \emph{quasi-open} when for any open set $U$ in $X$,
the image $f(U)$ has nonempty interior.
\edfn \mkp

We now formulate our main results. All topological spaces are assumed to be infinite, completely regular and Hausdorff from here on.

\begin{thmx}\label{theorem_A}
Let $X$ be a \v{C}ech-complete space, $M$ a metric space, $Y$ a metrizable separable space and $\Phi:X\rightarrow Y$ a continuous and quasi-open map.
If $G\subseteq C(X,M)$ is a $\frak{B}$-family such that each $g\in G$ factors through $\Phi$ (that is,
for each $g\in G$, there is a map $\widetilde{g}\in C(Y,M)$ satisfying $g(x)=(\widetilde{g}\circ \Phi)(x)$ for all $x\in X$),
then there is a nonempty compact subset $\Delta$ of $X$ and a countable infinite subset $L$ of $G$ such that $L$ is separated by $\Delta$.
As a consequence, if $M$ is a Banach space, $G$  contains a countably infinite $M$-interpolation set.
\end{thmx}

\brem From Theorem \ref{pol}, one can deduce the existence of an interpolation subset in a set $G$ of
real-valued continuous functions defined on a complete metric space $X$, when $\overline{G}^{\, X}$
contains a function that is not Baire one. The main difference in our approach is that
this property is isolated within the set $G$.
\erem

\begin{thmx}\label{theorem_B}
Let $X$ be a \v{C}ech-complete group and $K$ a compact group.
If $G\subseteq CHom(X,K)$ is not equicontinuous, then $G$ contains a countable subset $L$ such that 
$\overline{L}^{K^X}$ is canonically homeomorphic to $\beta L$, when $L$ is equipped with the discrete topology.
In case $K=\U(n)$, the unitary group of degree $n$, it follows that $L$ is an $\C^{n^2}$-interpolation set.
\end{thmx}
\mkp

A consequence of this result is a variation of a well-know Theorem by Corson and Glicksberg \cite{Corson1970} asserting that
if a subset $G$ of continuous homomorphisms defined on a hereditarily Baire group has a compact, metric closure,
then it is equicontinuous. In case $X$ is \v{C}ech-complete and $K$ is a compact group, these constraints can be relaxed considerably.

\begin{thmx}\label{theorem_C}
Let $X$ be a \v{C}ech-complete group, $K$ be a compact group and $G$ be an infinite subset of $CHom(X,K)$.
If for every countable subset
$L\subseteq G$ and compact separable subset $Y\subseteq X$ we have that either
$\overline{L}^{K^Y}$ has countable tightness or $\vert\overline{L}^{K^Y}\vert \leq \mathfrak{c} $, then $G$ is equicontinuous.
\end{thmx}

\bdfn
A Hausdorff topological space $X$ is a $k_{\omega}$\textit{-space} if there exists an ascending sequence of compact subsets
$K_1\subseteq K_2\subseteq \ldots\subseteq X$ such that $X=\bigcup\limits_{n<\omega}K_n$ and $U\subseteq X$ is open if and only if
$U\cap K_n$ is open in $K_n$ for each $n<\omega$ (i.e. $X=\lim\limits_{\rightarrow}K_n$) as a topological space.
A Hausdorff topological space $X$ is \textit{locally} $k_{\omega}$ if each point has an open neighbourhood which is a
$k_{\omega}$-space in the induced topology. It is clear that every $k_\omega$-space is a $k$-space (see \cite{Glockner2010}).
A $k_{\omega}$\textit{-group} (resp. \textit{locally} $k_{\omega}$\textit{-group}) is a topological group
where the underlying topological space is a $k_{\omega}$-space (resp.  locally $k_{\omega}$).

The class of abelian locally quasiconvex, locally $k_{\omega}$-groups includes, in addition to all locally compact abelian groups:
all free abelian groups on a compact space, indeed on any $k_\omega$ space; all dual groups of countable projective limits of metrizable
(more generally, \v{C}ech-complete) abelian groups; all dual groups of abelian pro-Lie groups defined by countable systems
\cite{Glockner2010,HofmannMorris_book2007}.
Moreover, this class is preserved by countable direct sums, closed subgroups, and finite products \cite{Glockner2010}.
\edfn

\begin{thmx}\label{theorem_D}
Let G be an abelian locally quasiconvex, locally $k_\omega$-group.
If $\lbrace g_n\rbrace_{n<\omega}$ is a sequence in $G$ that is not precompact in $G$,
then $\lbrace g_n\rbrace_{n<\omega}$ contains an $I_0$-set.
\end{thmx}

The Bohr compactification of a topological  group $G$,
can be defined as a pair $(bG,b)$ where $bG$ is
a compact Hausdorff group  and $b$ is a continuous homomorphism
 from $G$ onto a dense subgroup of $bG$ such that
every  continuous homomorphism $h\colon G\to K$ into a  compact
group $K$ extends to   a continuous homomorphism $h^{b}\colon bG
\to K$, making the following diagram
commutative:
\[
\xymatrix{ G \ar@{>}[rr]^{b} \ar[dr]^{h} & & bG \ar[dl]_{h^b} \\ & K &}
\]

The topology that $b$ induces on $G$ will be referred to as the
\textit{Bohr topology}. A topological group $G$ is said to be
\emph{maximally almost periodic} (MAP, for short) when the map $b$ is one-to-one,
which implies that the Bohr topology will be  Hausdorff.

The duality theory can be used to represent the Bohr compactification of an abelian group as a group of homomorphisms.
Indeed, if $G$ is an abelian
topological group and $\Gamma _{d}$ denotes its dual group
equipped with the discrete topology then $bG$ coincides with the dual group
of $\Gamma _{d}$. 

Given a topological group $G$, let $G^+$ denote
the algebraic group $G$ equipped with the Bohr topology.
Glicksberg \cite{Glicksberg1962} has shown that in a locally compact abelian (LCA, for short) group $G$,
every compact subset in $G^+$ is compact in $G$.
This result concerning LCA groups is one of the pivotal results of the
subject, often referred to as \emph{Glicksberg's theorem}.

Given a topological group $G$ and a property $\mathcal P$, we say after Trigos-Arrieta \cite{Trigos-Arrieta1991a} that
$G$ \emph{respects the property $\mathcal P$} when $G$ and $G^+$ have the same sets satisfying $\mathcal P$.
Taking this terminology, Glicksberg's theorem asserts that locally compact
Abelian groups respect compactness. Trigos-Arrieta considered some properties (pseudocompactness, countable
compactness, functional boundedness) obtaining that they are respected by locally compact
Abelian groups. Several authors have dealt with this question subsequently (cf. \cite{Ausenhofer2008,Banaszczyk1999,Hernandez2001a,Gabriyelyan2017}).
Glicksberg result was extended in a different direction by Comfort, Trigos-Arrieta and Wu \cite{Comfort1993}
by the following remarkable result.

Let $G$ be a LCA group and let $N$ be a closed metrizable subgroup of its Bohr compactification $bG$.
Denote by $\pi$ the canonical projection from $bG$ onto $bG/N$ and set
$b_N\defi \pi\circ b$ making the following diagram commutative:

\[
\xymatrix{ G \ar@{>}[rr]^{b} \ar[dr]^{b_N} & & bG \ar[dl]_{\pi} \\ & \frac{bG}{N} &}
\]

\bthm[Comfort, Trigos-Arrieta and Wu]\label{ctw}
Let $G$ be a LCA group and let $N$ be a closed metrizable subgroup of its Bohr compactification $bG$.
If $A$ is a subset of $G$, then  $A+(N\cap G)$ is compact in $G$ if and only if the set $b_N(A)$ is compact in $bG/N$.
\ethm

In the same paper, the following classes of topological grous is introduced:
A group $G$ \textit{strongly respects compactness}
if satisfies the thesis in Theorem \ref{ctw}.
The authors also propose the question of clarifying the relation between these two classes of groups
and furthermore the characterization of the groups that strongly respect compactness.
Using the techniques studied in this paper, we can prove that every abelian locally quasiconvex, locally $k_{\omega}$ group
respects any compact-like property $\mathcal P$ that implies functional boundedness and, furthermore,
strongly respects compactness, improving the results obtained by Gabriyelyan \cite{Gabriyelyan2017} for
locally $k_\omega$-groups. As a matter of fact, this result has been already applied to solve Question 4.1 in
\cite{Comfort1993} (see \cite{HdezTrigos2017}).

\bdfn
Let $X$ be a topological space. A subset $A$ of $X$ is \emph{functionally bounded} when every real-valued continuous
function defined on $X$ is bounded on $A$. We say that a topological property $\mathcal P$ on $X$ is
\emph{stronger than or equal to functional boundedness}
if for each $A\subseteq X$ that satisfies $\mathcal P$ ($A\in \mathcal P$ for short), it holds that $A$ is
functionally bounded.
\edfn

\begin{thmx}\label{theorem_E}
Let $G$ be an abelian, locally quasiconvex, locally $k_\omega$, group. Then the following holds:
\begin{enumerate}
\item[(i)] $G$ respects any compact-like property $\mathcal P$ stronger than or equal to functional boundedness.
\item[(ii)] $G$ strongly respects compactness.
\end{enumerate}
\end{thmx}

\section{Interpolation sets in topological spaces}
\bdfn\label{I-set}
Let $X$ and $M$ be a topological space and metric space, respectively, and let $C(X,M)$ denote
the space of continuous functions of $X$ into $M$.
Given a subset $L\subseteq C(X,M)$, we say that $K\subseteq X$ \emph{separates $L$}
if for every subset $A\subseteq L$ there are two closed subsets in $M$, say $D_{1}$ and $D_{2}$, and
$x_A \in K$ such that $\hbox{dist}(D_1,D_2)>0$, $\chi (x_A)\in D_{1}$ for all $\chi\in A$
and $\chi(x_A)\in D_{2}$ for all $\chi\in L\setminus A$.
\edfn
\mkp

In the sequel, we are going to apply the definition of $M$-interpolation set to subsets $L\subseteq C(X,M)\subseteq M^X$,
where $X$ and $M$ are a topological and a metric space, respectively. That is to say, we will look at $L$
as an Interpolation set for $C(M^X,M)$. First, we need a lemma, whose proof is known. However,
we include it here for the reader's sake. We refer to \cite{Engelking1989,Gillman1960,Prolla1977} for further information.

\blem\label{lem_i0}
Let $X$ and $M$ be a topological and a metric space, respectively, and let $L$ be a subset of $C(X,M)$
such that $\overline{L}^{M^X}$ is compact.
Consider the following properties:
\begin{enumerate}[(a)]
\item There is a nonempty subset $\Delta$ of $X$ such that $L$ is separated by $\Delta$.
\item Every two disjoint subsets of $L$ have disjoint closures in $M^X$.
\item $\overline{L}^{M^X}$ is canonically homeomorphic to $\beta L$ if $L$ is equipped with the discrete to\-po\-lo\-gy.
\item  $L$ is a Interpolation set for $C(M^X,M)$.
\end{enumerate}
Then $(a)\Rightarrow (b)\Leftrightarrow (c) \Leftarrow (d)$. If $M$ is a Banach space then the properties $(b)$, $(c)$ and $(d)$ are equivalent.
\elem
\bpf
That (b) implies (c) is folklore. It is also clear that $(d)$ implies $(c)$. For $(a)$ implies $(b)$, let $B_1$ and $B_2$ two disjoint subsets of $L$, which is separated by $\Delta$. Then, there are two closed sets $D_1$ and $D_2$ in $M$ and $x_0\in \Delta\subseteq X$ such that $d(D_1,D_2)\geq \epsilon_0$, for some $\epsilon_0>0$, $b_1(x_0)\in D_1$ for all $b_1\in B_1$ and $\gamma(x_0)\in D_2$ for all $\gamma\in L\setminus B_1$ (in particular for all $b_2\in B_2$).
Thus, $\overline{B_1}^{M^X}\cap \overline{B_2}^{M^X}=\emptyset$.
Finally, let us see that $(c)$ implies $(d)$, assuming that $M$ is a Banach space.

Let $f\in M^L$ with relatively compact range in $M$. By $(c)$, the map $f$ can be extended to a continuous map defined on $\overline{L}^{M^X}$.
Therefore, there is a continuous function ${\overline f \colon \overline{L}^{M^X} \to M}$ such that $\overline f_{|L}=f$.
Now, applying \cite[Cor. 3.5]{Prolla1977} (cf. \cite[Th. 9]{Hernandez-Munoz1994}), it follows that there is a continuous map
$\widetilde{f}\colon {M^X}\to M$ that extends $\overline f$. Hence $\widetilde{f}$ is the required extension of $f$ to $M^X$.
\epf

\bdfn\label{def_totally discontinuous}
Let $X$ and $M$ be a topological space and a metric space (respectively) and let $f\in M^X$.
We say that $f$ is \emph{totally discontinuous} if there are two subsets $N_0$ and $N_1$ in $M$
and two dense subsets $A_0$ and $A_1$ in $X$ such that $d(N_0,N_1)>0$ and $f(A_j)\subseteq N_j$ for $j=0,1$.
\edfn

We may assume that $N_0$ and $N_1$ are open sets because, otherwise, we would replace them by
$B(N_i,s/3)\defi\lbrace m\in M: d(m,N_i)<s/3\rbrace$, where $s=d(N_0,N_1)$ and $i=0,1$.

\bdfn\label{def_cech_comple}
A topological space X is said to be \emph{\v{C}ech-complete} if it is a $G_{\delta}$-subset of its Stone-\v{C}ech compatification.
The family of \v{C}ech-complete spaces includes all complete metric spaces and all
locally compact spaces.
\edfn

\blem\label{lem_resultado_5}
Let $X$ and $M$ be a \v{C}ech-complete space and a metric space, respectively. If $G$ a subset of $C(X,M)$ where each element has relatively compact range in $M$ such that $\overline{G}^{\, M^X}$ contains a totally discontinuous function $f$,
then there is a nonempty compact subset $\Delta$ of $X$ and a countable infinite subset $L$ of $G$, which is separated by $\Delta$.
Furthermore, by Lemma \ref{lem_i0}, if $M$ is a Banach space, it follows that $L$ is a $M$-interpolation set.
\elem
\bpf
Since $X$ is \v{C}ech-complete, it is a $G_\delta$-subset of its Stone-Cech compatification $\beta X$.
Set $X=\bigcap\limits_{n=0}^{\infty}W_n$, where $W_n$ is a dense open subset of $\beta X$
for each $n<\omega$ and $W_{s}\subseteq W_{r}$ if $r<s$. In the sequel,
given a map $g\in C(X,M)$ with relatively compact range in $M$, we denote by $g^{\beta}$ its continuous extension to $\beta X$.

Set $N_0, N_1, A_0, A_1$ as in Definition \ref{def_totally discontinuous}, where we assume that $N_0$ and $N_1$ are open wlog.
By induction on $n=\vert t \vert$, $t\in 2^{(\omega)}$ (i.e. the set of finite sequences of $0$'s and $1$'s),
we define a family $\lbrace U_t:t\in 2^{(\omega)}\rbrace$ of non-empty open subsets in $\beta X$
and a sequence of functions $\lbrace h_n:n<\omega\rbrace\subseteq G$, satisfying the following conditions for all $t\in 2^{(\omega)}$:

  \begin{enumerate}[(i)]
  \item $U_{\emptyset}\subseteq \overline{U_{\emptyset}}^{\beta X}\subseteq W_0$;
  \item $U_{ti}\subseteq \overline{U_{ti}}^{\beta X}\subseteq W_{\vert t\vert+1}\cap U_t$ for $i=0,1$;
  \item $U_{t0}\cap U_{t1}=\emptyset$;
  \item $h^{\beta}_{\vert t\vert}(U_{tj})\subseteq {N}_j$ for $j=0,1$;
  \item if $s<\vert t \vert$, then $diam(h^{\beta}_s(\overline{U_{tj}}^{\beta X}))<\frac{1}{\vert t\vert}$ for $j=0,1$.
  \end{enumerate}

\textit{Construction:} If $n=0$, by regularity we can find $U_{\emptyset}$ a nonempty open set in $\beta X$ such that
$U_{\emptyset}\subseteq \overline{U_{\emptyset}}^{\beta X}\subseteq W_0$. For $n\geq 0$,
suppose $\lbrace U_t:\vert t\vert \leq n\rbrace$ and $\lbrace h^{\beta}_{\vert t\vert}:\vert t\vert <n\rbrace$ have been defined satisfying $(i)-(v)$.
Fix $t\in 2^{(\omega)}$ with $\vert t\vert =n$. Since $U_t$ is open in $\beta X$, then $V_t\defi U_t\cap X\neq \emptyset$ is open in $X$.
We can find $a_t,b_t\in V_t$ such that $f(a_t)\in {N}_0$ and $f(b_t)\in {N}_1$ because $V_t$ is a relatively open subset of $X$
and the sets $A_0$ and $A_1$ are dense in $X$.
Since $f\in \overline{G}^{M^X}$, there is $h_n\in G$  such that $h_n(a_t)\in {N}_0$ and $h_n(b_t)\in {N}_1$.

Let $h_n^{\beta}$ be the continuous extension of $h_n$, then we can select two open disjoint neighborhoods in $\beta X$,
$O_{t0}$ and $O_{t1}$ of $a_t$ and $b_t$,
respectively, satisfying:

 \begin{enumerate}[(1)]
\item $\overline{O_{t0}}^{\beta X}\cup \overline{O_{t1}}^{\beta X}\subseteq U_t$;
\item $\overline{O_{t0}}^{\beta X}\cap \overline{O_{t1}}^{\beta X}=\emptyset$;
\item $diam(h^{\beta}_n(O_{tj}))<\frac{1}{\vert t \vert}$;
\item $ h^{\beta}_n(O_{tj})\subseteq {N_j},\text{ for }j=0,1$.
\end{enumerate}

Since $W_{\vert t \vert +1}$ is dense in $\beta X$, then $W_{\vert t \vert +1}\cap O_{t0}$ and $W_{\vert t \vert +1}\cap O_{t1}$
are two nonempty open sets. By regularity, there exist
two non empty open sets $U_{t0}$ and $U_{t1}$ such that $U_{t0}\subseteq \overline{U_{t0}}^{\beta X}\subseteq W_{\vert t\vert +1}\cap O_{t0}$ and
$U_{t1}\subseteq \overline{U_{t1}}^{\beta X}\subseteq W_{\vert t\vert +1}\cap O_{t1}$, respectively.
Therefore, $U_{t0}$ and $U_{t1}$ satisfies the conditions $(ii),(iii)$ and $(iv)$. Moreover,
using a continuity argument, we can adjust the two open sets to satisfy $(v)$.

Set $\Delta\defi \bigcap\limits_{n<\omega}\bigcup\limits_{\vert t \vert=n}\overline{U_t}^{\beta X}$, which is a closed subset of $\beta X$ and,
as a consequence, $\Delta$ is compact. On the other hand, we also have
$\Delta=\bigcup\limits_{\sigma\in 2^{\omega}}\bigcap\limits_{n<\omega}\overline{U_{\sigma\vert n}}^{\beta X}$.
For each $\sigma\in 2^{\omega}$, $\bigcap\limits_{n<\omega}\overline{U_{\sigma\vert n}}^{\beta X}\neq \emptyset$ by compactness of $\beta X$.
So $\Delta\neq \emptyset$. By construction we have that $\Delta\subseteq \bigcap\limits_{n=0}^{\infty}W_n=X$.
Consequently $\Delta$ is contained in $X$.

Define $\varphi:\Delta\rightarrow 2^{\omega}$ by $\varphi^{-1}(\sigma)=\bigcap\limits_{n<\omega}\overline{U_{\sigma\vert n}}^{\beta X}$.
Clearly $\varphi$ is an onto and continuous map. For each $t\in 2^{(\omega)}$ and $\sigma\in 2^{\omega}$,
$h_{|t|}(\varphi^{-1}(\sigma))$ is a singleton by $(v)$. Therefore, $h_{|t|}$ lifts to a continuous function $h_{|t|}^*$ on $2^{\omega}$ such that
$h_{|t|}(x)=h_{|t|}^*(\varphi(x))$ for all $x\in \Delta$.

Let us see that  $\lbrace h_n\rbrace_{n<\omega}$ is separated by $\Delta$. Indeed, for any arbitrary subset $S\subseteq \omega$, it suffices to select
$\sigma\in 2^{\omega}$ such that $\sigma(0)=0$ and $\sigma(n+1)=1$ if $n\in S$ or $\sigma(n+1)=0$ if $n\not\in S$. By construction,
if we take any element $z\in \bigcap\limits_{n<\omega}\overline{U_{\sigma\vert n}}^{\beta X}\subseteq \Delta$,
then  $h_n(z)\in {N}_1$ for every $n\in S$ and $h_n(z)\in {N}_0$ for every $n\not\in S$.
Finally, in case $M$ is a Banach space, Lemma \ref{lem_i0} implies that $L=\lbrace h_n\rbrace_{n<\omega}$ is a $M$-interpolation set.
\epf
\mkp

We need the following compact space $\mathcal{K}$ that is defined as in \cite{Christensen1981}.

\bdfn
Let $(M,d)$ be a metric space that we always assume equipped with a bounded metric. We set
$$ \mathcal{K}\defi\lbrace \alpha:M\rightarrow [-1,1]:\vert \alpha(m_1)-\alpha(m_2)\vert\leq d(m_1,m_2),\quad\forall m_1,m_2\in M\rbrace.$$

\noindent Being pointwise closed and equicontinuous by definition, it follows that $\mathcal{K}$ is a compact and metrizable subspace of $C(M,\mathbb R)$
equipped with the supremum norm.
For each $m_0\in M$, we set $\alpha_{m_0}\in \mathbb R^M$ by $\alpha_{m_0}(m)\defi d(m,{m_0})$ for all $m\in M$.
It is easy to check that $\alpha_{m_0}\in \mathcal{K}$.
Given any element $f\in M^X$, we associate a map
$\check{f}\in \mathbb R^{X\times \mathcal{K}}$ defined by $$\check{f}(x,\alpha)=\alpha(f(x))\ \hbox{for all}\ (x,\alpha)\in X\times \mathcal{K}.$$
In like manner, given any subset $G$ of $M^X$ we set $\check{G}\defi\lbrace \check{f}:f\in G\rbrace$.
\edfn

We are now in position of proving  Theorem \ref{theorem_A}.

\begin{proof}[\textbf{Proof of Theorem \ref{theorem_A}}]
{We may assume wlog that the map $\Phi$ is surjective because otherwise we would deal with the separable and metrizable space $\Phi(X)$.}
Due to the fact that \v{C}ech-completeness is hereditary for closed subsets, we may assume, from here on, that
$X=\overline{V}$ wlog; where $V$ is a nonempty open subset satisfying the following property: there is some fixed $\epsilon>0$ such that
for every finite collection $\lbrace U_1,\ldots,U_n\rbrace$ of non-empty open subsets contained in $V$,
there is some element $g\in G$ with $diam(g(U_{j}))\geq\epsilon$ for all
$j\in \lbrace 1,\ldots,n\rbrace$.

Let $\lbrace \widetilde{V_k}\rbrace_{k<\omega}$ be an arbitrary countable open basis in $Y$. We set
$V_k\defi \Phi^{-1}(\widetilde{V_k})$ and pick and arbitrary point $x_k\in V_k$ for each $k<\omega$.

Since $X$ is \v{C}ech-complete, there exists a sequence $\lbrace\mathcal{A}_i\rbrace_{i<\omega}$ of open coverings of $X$,
such that, if a family $\mathcal{F}$ of closed subsets has the finite intersection property, and if for each $i<\omega$ there is an element of
$\mathcal{F}$ such that is contained in a member of $\mathcal{A}_i$, then $\bigcap\mathcal{F}\neq \emptyset$ \cite[Theorem 3.9.2]{Engelking1989}.
In order to simplify the notation below, we say that a set of $X$ is $\mathcal{A}_i$-\emph{small} if it is contained in a member of $\mathcal{A}_i$.

Using an inductive argument, for every integer $n<\omega$, we find $f_n\in G$, $\alpha_n\in\mathcal{K}$ and
a finite collection $\lbrace U_{n,k}\rbrace_{1\leq k\leq n}$
of nonempty open sets in $X$ satisfying the following conditions (for each $n<\omega$ and each $k=1,\ldots,n$):
  \begin{enumerate}[(i)]
  \item $U_{n,k}\subseteq V_k$;
  \item $diam(f_n(U_{n,k}))\leq \frac{1}{n}$;
  \item $\overline{U_{n+1,k}}\subseteq U_{n,k}$;
  \item $d( f_n(x),f_n(x_k))\geq \frac{\epsilon}{3}$, for all $x\in U_{n,k}$;
  \item $U_{n,k}$ is $\mathcal{A}_j$-small, for $1\leq j\leq n$;
  \item $\vert \alpha_n(f_n(x))-\alpha_n(f_n(x_k)\vert \geq \frac{\epsilon}{3}$, for all $x\in U_{n,k}$.
  \end{enumerate}

\textit{Construction:} If $n=1$, since $V_1$ is an open subset in $X$ there exists $f_1\in G$ such that $diam(f_1(V_1))\geq\epsilon$.
By the continuity of $f_1$, it follows that there exists a nonempty open subset $W_{1,1}$  such that:
  \begin{enumerate}[(a)]
  \item $W_{1,1}\subseteq V_1$
  \item $d( f_1(x),f_1(x_1))\geq \frac{\epsilon}{3}$, for all $x\in W_{1,1}$
  \end{enumerate}
Let $\alpha_1\defi \alpha_{f_1(x_1)}\in \mathcal{K}$. Note that $\vert \alpha_1(f_1(x))-\alpha_1(f_1(x_1))\vert\geq \frac{\epsilon}{3}$, for all $x\in W_{1,1}$.

Now, we take the open covering $\mathcal{A}_1$ of $X$. Then, there is $A\in \mathcal{A}_1$ such that
$A\cap W_{1,1}$ is not empty. By regularity, we can find a nonempty open subset $U_{1,1}$ such that
$U_{1,1}\subseteq\overline{U_{1,1}}\subseteq A\cap W_{1,1}\subseteq V_1$ and $diam(f_1(U_{1,1})\leq 1$.

Assume now that $f_n$, $\alpha_n$ and $\lbrace U_{n,k}\rbrace_{1\leq k\leq n}$ have been obtained, with $n\geq 1$.
By hypothesis, there exists $f_{n+1}\in G$ such that $diam(f_{n+1}(U_{n,k}))\geq\epsilon$
for all $k\in \lbrace 1,\ldots, n\rbrace$ and $diam(f_{n+1}(V'_{n+1}))\geq\epsilon$, where $x_{n+1}\in V_{n+1}'\subseteq V_{n+1}$ and $V_{n+1}'$ is $\mathcal{A}_j$-small for $1\leq j\leq n$.

By the continuity of $f_{n+1}$, we can find nonempty open subsets $\lbrace W_{n+1,k}\rbrace_{1\leq k\leq n+1}$ 
satisfying:
  \begin{enumerate}
  \item $W_{n+1,k}\subseteq U_{n,k}$, for all $1\leq k\leq n$;
  \item $W_{n+1,n+1}\subseteq V'_{n+1}$ (therefore $W_{n+1,n+1}$ is $\mathcal{A}_j$-small for $1\leq j\leq n$);
  \item $diam(f_{n+1}(W_{n+1,k}))\leq \frac{1}{n+1}$, for all $1\leq k\leq n+1$;
  \item $d( f_{n+1}(x),f_{n+1}(x_k))\geq \frac{\epsilon}{3}$, for all $x\in W_{n+1,k}$ and $1\leq k\leq n+1$.
  \end{enumerate}

Set $\alpha_{n+1}\in [-1,1]^M$ defined by $$\alpha_{n+1}(m)\defi \min\limits_{1\leq k\leq n+1} d(m,f_{n+1}(x_k)) \text{ for all } m\in M.$$
We claim that $\alpha_{n+1}\in \mathcal{K}$. Indeed, if $m_1,m_2\in M$,
then $$\vert \alpha_{n+1}(m_1)-\alpha_{n+1}(m_2)\vert= \vert \min\limits_{1\leq k\leq n+1} d(m_1,f_{n+1}(x_k))-\min\limits_{1\leq k\leq n+1} d(m_2,f_{n+1}(x_k))\vert.$$
Assume wlog that $$\min\limits_{1\leq k\leq n+1} d(m_1,f_{n+1}(x_k))\geq\min\limits_{1\leq k\leq n+1} d(m_2,f_{n+1}(x_k))$$
and choose $k_0\in\lbrace 1,\ldots,n+1\rbrace$ such that $$\min\limits_{1\leq k\leq n+1} d(m_2,f_{n+1}(x_k))=d(m_2,f_{n+1}(x_{k_0})).$$ Then,
\begin{equation*}
\begin{array}{c}
\vert \alpha_{n+1}(m_1)-\alpha_{n+1}(m_2)\vert= \min\limits_{1\leq k\leq n+1} d(m_1,f_{n+1}(x_k))-d(m_2,f_{n+1}(x_{k_0}))\leq\\
\\
d(m_1,f_{n+1}(x_{k_0}))-d(m_2,f_{n+1}(x_{k_0}))\leq  d(m_1,m_2).
\end{array}
\end{equation*}

On the other hand, for all $x\in W_{n+1,k'}$ and $1\leq k'\leq n+1$:
 $$\vert \alpha_{n+1}(f_{n+1}(x))-\alpha_{n+1}(f_{n+1}(x_{k'}))\vert=\vert \alpha_{n+1}(f_{n+1}(x))\vert=\min\limits_{1\leq k\leq n+1} d(f_{n+1}(x),f_{n+1}(x_k))\geq \frac{\epsilon}{3},$$

Take the open covering $\mathcal{A}_{n+1}$ of $X$. Then, for each $k\in \lbrace 1,\ldots, n+1\rbrace$ there is $A_k\in \mathcal{A}_{n+1}$
such that $A_{k}\cap W_{n+1,k}$ is a nonempty open subset of $X$.
By regularity we can find an open set $U_{n+1,k}$ such that:

\begin{enumerate}[$\bullet$]
 \item $U_{n+1,k}\subseteq\overline{U_{n+1,k}}\subseteq A_k\cap W_{n+1,k}\subseteq U_{n,k},\quad\text{ if }1\leq k\leq n$;
 \item $U_{n+1,n+1}\subseteq\overline{U_{n+1,n+1}}\subseteq A_{n+1}\cap W_{n+1,n+1}\subseteq V_{n+1},\quad\text{ if }k=n+1$.
\end{enumerate}
This completes the construction.

Now, for each $k<\omega$, the intersection $\bigcap\limits_{n=k}^{\infty}U_{n,k}$ is nonempty by \v{C}ech-completeness.
Therefore, we can fix a point $z_{k}\in \bigcap\limits_{n=k}^{\infty}U_{n,k}$ for all $k<\omega$.
Note that $\Phi(x_k)\in \widetilde{V_k}$ and $\Phi(z_k)\in \widetilde{V_k}$ for all $k\in \omega$.

Take an element $(f,\alpha)\in \overline{\lbrace (f_n,\alpha_n)\rbrace_{n<\omega}}^{(\overline{G}^{M^X}\times \mathcal{K})}$.\\
By $(vi)$ we have:
$$\vert\alpha_n\circ \widetilde{f_n}(\Phi(z_k))-\alpha_n\circ \widetilde{f_n}(\Phi(x_k))\vert=\vert\alpha_n\circ f_n(z_k)-\alpha_n\circ
f_n(x_k)\vert\geq\frac{\epsilon}{3},\quad\forall n\geq k.$$

Therefore, $osc(\alpha_n\circ \widetilde{f_n},\widetilde{V_k})\geq\frac{\epsilon}{3}$ for all $n\geq k$.
As a consequence, we also have $osc(\alpha\circ \widetilde{f},\widetilde{V_k})\geq\frac{\epsilon}{3}$ for all $k< \omega$.

Let $\lbrace r_m, \delta_m\rbrace_{m<\omega}$ be an enumeration of all pairs of rational numbers $(r,\delta)$ with $\delta>0$.
For each $m<\omega$, define
$$\widetilde{F_m}=\lbrace y\in Y:\inf(\alpha\circ \widetilde{f})(U)<r_m,\quad \sup(\alpha\circ \widetilde{f})(U)\geq r_m+\delta_m,
\forall \text{nbd }U\text{ of }y\rbrace.$$
It is easily seen that $\widetilde{F_m}$ is closed and, consequently, $F_m\defi\Phi^{-1}(\widetilde{F_m})$ is closed in $X$.\\
Observe that, since $\lbrace \widetilde{V_k}\rbrace_{k<\omega}$ is an open basis in $Y$, it follows that $Y=\bigcup\limits_{m<\omega}\widetilde{F_m}$ and,
hence $X=\bigcup\limits_{m<\omega}F_m$. Being $X$ \v{C}ech-complete, it is a Baire space.
Therefore, there is some $m_0<\omega$ such that $F_{m_0}$ has nonempty interior $U$ in $X$. Since $\Phi$ is a quasi-open map, we have that
$\overline{\Phi(U)}$ has nonempty interior $\widetilde{U}$ included in $\widetilde{F}_{m_0}$. It follows that
$\inf \left( \alpha\circ \widetilde{f}(\widetilde{U})\right) <r_{m_0}$ and $\quad \sup \left( \alpha\circ \widetilde{f}(\widetilde{U})\right) \geq r_{m_0}+\delta_{m_0}$.
Set $U_0=\Phi^{-1}(\widetilde{U})\subseteq U$
we have that $\inf \left( \alpha\circ f(U_0)\right) <r_{m_0}$ and $\sup \left( \alpha\circ f(U_0)\right) \geq r_{m_0}+\delta_{m_0}$.

Set $F=\overline{U_0}$, $r\defi r_{m_0}$ and $\delta\defi \delta_{m_0}$ and consider the following sets:
\[
A_0=\lbrace x\in F: \alpha\circ f(x)<r\rbrace=\lbrace x\in F: \alpha\circ f(x)\in I_0\rbrace\]
\[
A_1=\lbrace x\in F: \alpha\circ f(x)\geq r+\delta\rbrace=\lbrace x\in F: \alpha\circ f(x)\in I_1\rbrace
\]

where $I_0=[-1,r)$ and $I_1=(r+\delta,1]$. Note that $A_0$ and $A_1$ are dense subsets in $F$.
Define ${N_0}\defi\alpha^{-1}(I_0)$ and ${N_1}\defi\alpha^{-1}(I_1)$, which are disjoints. Moreover, since $\alpha\in\mathcal K$,
it follows that $d({N_0},{N_1})\geq \delta$ and $f(A_j)\subseteq N_j$ for $j=0,1$. Therefore $f$ is totally discontinuous on $F$.
It now suffices to apply Lemma \ref{lem_resultado_5}
\end{proof}
\bkp

\brem\label{nota_5}
Note that the result remains valid if we assume that for each residual subset $R$ of $X$ there is a separable
metrizable space $Y$ and a continuous and quasi-open map $\Phi:R\rightarrow Y$ such that for all $g\in G$ there is a
$\widetilde{g}\in C(Y,M)$ satisfying $g(x)=(\widetilde{g}\circ \Phi)(x)$ for all $x\in R$.
\erem

\bcor\label{cor_5_polaco}
Let $X$ be a polish space, let $(M,d)$ be a metric space and let $G\subseteq C(X,M)$ be a $\frak{B}$-family. Then there is a nonempty compact subset $\Delta$ of $X$ and a countable subset $L$ of $G$ such that $L$ is separated by $\Delta$. As a consequence, if $M$ is a Banach space, it follows that $L$ is a $M$-interpolation set.
\ecor
\mkp

\section{Interpolation sets in topological groups}

In this section, we apply the results obtained previously in the setting of topological groups.
Our first result clarifies the relevance of the notion of $\frak B$-family when we deal with topological groups.
From here on, we assume, wlog, that every metrizable topological group $M$ is equipped with a left-invariant
metric. Furthermore, if $M$ is in addition compact, then we assume that $M$ is equipped with a bi-invariant metric.
From here on, if $X$ and $Y$ are topological groups, we let $Hom(X,M)$ (resp. $CHom(X,M)$) denote
the set of all homomorphisms (resp. continuous homomorphisms) of $X$ into $M$.

\blem\label{resultado_6}
Let $X$ be a topological group, $M$ a metric topological group and $G\subseteq CHom(X,M)$ such that $\overline{G}^{M^X}$ is compact. Then $G$ is a $\frak B$-family if and only if it is not equicontinuous.
\elem
\begin{proof}
It is clear that, if $G$ is a $\frak B$-family, then it may not be equicontinuous. So, assume that $G$ is not a $\frak B$-family.
Taking $V=X$ and $\epsilon>0$ arbitrary, there exists a finite family $\lbrace U_1,\ldots,U_n\rbrace$ open subsets in $X$
(wlog, we assume that $U_j=x_jV_j$, where $V_j$ is a neighbourhood of the neutral element) such that
for every $g\in G$ there is $V_j$, with $1\leq j\leq n$, satisfying that $diam(g(x_jV_j))<\epsilon$.
Now, since $g$ is a group homomorphism and $d$ is left-invariant, it follows that $diam(g(V_j))<\epsilon$ as well.
Set $V_0=V_1\cap\ldots\cap V_n$, then $diam(g(xV_0))<\epsilon$ for all $g\in G$ and $x\in X$. Consequently $G$ is equicontinuous.
\end{proof}
\mkp

The next result is a direct consequence of Lemma \ref{resultado_6}, Theorem \ref{theorem_A} and Lemma \ref{lem_i0}.
Previously, we need the following definition.
Recall that $\U(n)$ denotes the unitary group of degree $n$.



\bcor\label{resultado_5_COR_GT}
Let $X$ be a compact group, $M$ a metric topological group and $G\subseteq CHom(X,M)$ such that $\overline{G}^{M^X}$ is compact.
If $G$ is not equicontinuous, then there is a nonempty compact subset $\Delta$ of $X$
and a countable infinite subset $L$ of $G$ that is separated by $\Delta$.
As a consequence, if $M$ is a Banach space, it follows that $L$ is a $M$-interpolation set.
In particular, if $M=\U(n)$ then $G$ contains an Interpolation set for $C(Hom(G,\U(n)),\C^{n^2})$.
\ecor
\begin{proof}
By Troallic \cite{Troallic1996}, we may assume wlog that $G$ is countable. By Lemma \ref{resultado_6}, $G$ is a $\frak{B}$-family.
Define an equivalence relation on $X$ by $x\sim y$ if and only if $g(x)=g(y)$ for all $g\in G$.
Since $G$ is countable and consists of group homomorphisms, it follows that the quotient space
$\widetilde{X}=X/{\sim}$ is a compact metrizable group. Therefore, if $p:X\rightarrow \widetilde{X}$ denotes the canonical quotient map,
each $g\in G$ factors through a map $\widetilde{g}$ defined on $CHom(\widetilde{X},M)$; that is $\widetilde{g}(p(x))\defi g(x)$ for any $x\in X$.
Since every quotient group homomorphism is automatically open, Theorem \ref{theorem_A} implies that
there is a nonempty subset $\Delta$ of $X$ and a subset $L$ of $G$ such that $L$ is separated by $\Delta$.
In case $M$ is a Banach space, applying Lemma \ref{lem_i0}, we obtain that $L$ is a $M$-interpolation set.
\end{proof}

Next result is folklore but we include its proof for the sake of completeness.

\blem\label{resultado_6_1}
Let $X$ be a topological group, $M$ a metric topological group, $G\subseteq C(X,M)$ and $h\in C(X,M)$.
Set $Gh\defi \lbrace gh:g\in G\rbrace$. Then {$G$ is equicontinuous on $X$ if and only if $Gh$ is equicontinuous on $X$.}
\elem
\bpf
It suffices to prove that $Gh$ is equicontinuous if $G$ is equicontinuous.
Let $x_0$ be an arbitrary but fixed point in $X$. Since right translations are continuous mappings on a topological group, and
$G$ (resp. $h$) is equicontinuous (resp. continuous) on $X$, given $\gep >0$,
there is a neighbourhood $U$ of $x_0$ such that $d(g(x_0)h(x_0),g(x)h(x_0))<\gep/2$
and $d(h(x_0),h(x))<\gep/2$ for all $x\in U$
and all $g\in G$.  Thus, applying the left invariance of the group metric, we obtain
$$d(g(x_0)h(x_0),g(x)h(x))\leq d(g(x_0)h(x_0),g(x)h(x_0))+d(g(x)h(x_0),g(x)h(x))<\frac{\gep}{2}+\frac{\gep}{2}=\gep.$$
for all $x\in U$, which completes the proof.
\epf
\mkp

With the hypothesis of the previous lemma, if $g\in CHom(X,M)$, the symbol $g^{-1}$ denotes the map defined by
$g^{-1}(x)=g(x)^{-1}=g(x^{-1})$ for all $x\in X$. Combining Lemmata \ref{resultado_6} and \ref{resultado_6_1}, we obtain:

\bcor\label{resultado_6_2}
Let $X$ be a topological group, $M$ be a topological group with a bi-invariant metric, $G\subseteq CHom(X,M)$ such that $\overline{G}^{M^X}$ is compact and $g_0\in G$.
Then $Gg_0^{-1}$ is a $\frak B$-family if and only if it is not equicontinuous.
\ecor
\begin{proof}
It suffices to see that $Gg_0^{-1}$ is equicontinuous if $Gg_0^{-1}$ is not a $\frak B$-family. Reasoning as in Lemma \ref{resultado_6_1}, let $V=X$ and $\epsilon>0$, then there are $\lbrace U_1,\ldots,U_n\rbrace$ open subsets of $X$
such that for all $g\in G$ there is $j\in \lbrace 1,\ldots,n\rbrace$ with $diam(gg_0^{-1}(U_j))<\epsilon$.
We can assume that $U_j=x_jV_j$ wlog, where $V_j$ is a neighbourhood of the identity element of $X$,
for all $1\leq j\leq n$. Take $W\defi \bigcap\limits_{1\leq j\leq n} V_j$ and an arbitrary element $x_0\in X$.

Given $g\in G$, there is $j\in \lbrace 1,\ldots,n\rbrace$ such that
\begin{eqnarray*}
\epsilon > diam( gg_0^{-1}(U_j))&=&diam( gg_0^{-1}(x_jV_j))=\sup\limits_{x,y\in V_j}d(gg_0^{-1}(x_jx),gg_0^{-1}(x_jy))\\& = &
\sup\limits_{x,y\in V_j} d(g(x_j)\left[gg_0^{-1}(x) \right]g_0^{-1}(x_j),g(x_j)\left[gg_0^{-1}(y) \right]g_0^{-1}(x_j))\\& = &
\sup\limits_{x,y\in V_j}d(\left[gg_0^{-1}(x) \right],\left[gg_0^{-1}(y) \right])=diam( gg_0^{-1}(V_j))\\&
\geq & diam( gg_0^{-1}(W))=diam( gg_0^{-1}(x_0W)).
\end{eqnarray*}
\end{proof}

We can now prove Theorem \ref{theorem_B}.

\begin{proof}[\textbf{Proof of Theorem \ref{theorem_B}}]
Since $K$ is compact, there is a representation $\pi:K\rightarrow \U(n)$ such that $\lbrace \pi\circ g:g\in G \rbrace$ is not equicontinuous.
Therefore, we assume that $K=\U(n)$ wlog.

 Applying \cite[Cor. 2.4]{Fer_Her_Tar2017} (cf. \cite[Cor. 3.2]{Troallic1996}), since $G\subseteq CHom(X,\U(n))$ is not equicontinuous, there exists
a separable compact subset $F$ of $X$ and a countable subset $L\subseteq G$ such that $L\vert_F$ is not equicontinuous. Set $H$ as the smallest closed
 subgroup generated by $F$, it follows that $H\leq X$ is closed and separable and $L\subseteq G$ countable such that $L\vert_H$ is not equicontinuous.
 So we can assume wlog that $X$ is separable and $G$ is countable. On the other hand, by \v{C}ech-completeness of $X$, there must be a
 compact subgroup $C$ of $X$ such that $X/C$ is separable, complete and metrizable \cite{Brown1972}, thereby , a Polish space.

Let $G\vert_C\defi \lbrace g\vert_C:g\in G\rbrace\subseteq CHom(C,\U(n))$. We have two possible cases:
\begin{enumerate}[(1)]
\item $G\vert_C$ contains infinitely many elements that are pairwise inequivalent
(recall that $\gamma_1,\gamma_2\in Hom(C,\U(n))$ are equivalent $(\gamma_1\sim\gamma_2)$ if there exists $U\in \U(n)$ such that $\gamma_1=U^{-1}\gamma_2 U$).
\item $G\vert_C$ only contains a finite subset of elements that are pairwise inequivalent.
\end{enumerate}

\begin{enumerate}[-]
\item Case $(1)$: We may suppose wlog that all elements of $G\vert_C$ are pairwise inequivalent.
By \cite[Th. 1]{Chu1966}, it follows that $G\vert_C$ is discrete as a subset of $CHom(C,\U(n))$ in the compact open topology on $C$,
which implies that $G\vert_C$ may not be equicontinuous on $C$. Applying Corollary \ref{resultado_5_COR_GT}, there is a nonempty subset $\Delta$
of $C$ and a countable subset $L$ of $G$ such that $L$ is separated by $\Delta$. Thus, by Lemma \ref{lem_i0}, $\overline L^{\U(n)^X}$
is canonically homeomorphic to $\beta L$ (where $L$ is equipped with the discrete topology) and we are done.

\item Case $(2)$: Set $H\defi \lbrace \phi_1,\ldots, \phi_m\rbrace\subseteq G$ such that every $g\in G$ is equivalent to an element in $H$ when they are restricted to $C$.
If we define $G_i=\lbrace g\in G:g\vert_C\sim \phi_{i}\vert_C\rbrace$, then $G=G_1\bigcup\ldots\bigcup G_m$. Since $G$ is not equicontinuous, there is
    $i\in\{1,\dots m\}$ such that $G_i$ is not equicontinuous. So, we may assume wlog that there is $g_0\in G$ such that
$g\vert_C\sim g_0\vert_C$ for all $g\in G$. Therefore, for each $g\in G$, there is $U_g\in \U(n)$ with $(U^{-1}_ggU_g)\vert_C=g_0\vert_C$.
Denote by $\widetilde{g}$ the map $U^{-1}_ggU_g$ and set $\widetilde{G}\defi \lbrace U^{-1}_ggU_g:g\in G \rbrace$, which is a subset of $CHom(X,\U(n))$.
It is easily seen that $\widetilde{G}$ is not equicontinuous on $X$.
(Indeed, assume that $\widetilde{G}$ were equicontinuous and let $W$ be an open neighbourhood of the identity matrix $I_{n}$ in $\U(n)$.
By \cite[Corollary 1.12]{Hofmann2006}, there would exist an open neighbourhood $V$ of $e_{X}$  such that
$\widetilde{g}(V)\subseteq \bigcap\limits_{U\in \U(n)}U^{-1}WU$ for all $\widetilde{g}\in\widetilde{G}$. Therefore, we would have
$g(v)=U_g \widetilde{g}(v)U_g^{-1}\in W$ for all $v\in V$. This would imply that $G$ is equicontinuous, which is a contradiction).
Hence,  $\widetilde{G}g_0^{-1}$ is  a $\frak B$-family on $X$ by Lemmas \ref{resultado_6_1} and \ref{resultado_6_2}.

Let $\pi_C:X\rightarrow X/C$ the canonical quotient map, which is open and continuous. Since $X/C$ is Polish and each $\widetilde{g}g_0^{-1}$ factors through $X/C$,
we apply Theorem \ref{theorem_A} and Lemma \ref{lem_i0} in order to obtain $\Delta \subseteq X$ and $\widetilde{L}\subseteq \widetilde G$ such that $$\overline{\widetilde{L}}^{\U(n)^{\Delta}}\simeq\overline{\widetilde{L}}^{\U(n)^{\Delta}}g_0^{-1}\simeq \beta \widetilde{L}.$$
Set $L\defi \lbrace g:\widetilde{g}\in \widetilde{L}\rbrace\subseteq G$ and consider the map
\begin{eqnarray*}
\psi:(\widetilde{L},t_p(\Delta)) &\longrightarrow & (L,t_p(\Delta))\\
  U_g^{-1}gU_g &\longmapsto &  g
\end{eqnarray*}
The map $\psi$ is continuous because $\widetilde{L}$ is {discrete}. Moreover, using that $\overline{\widetilde{L}}^{\U(n)^X}$ is canonically
 homeomorphic to $\beta L$ ($L$ with the discrete topology), there is a continuous extension map
 $$\overline{\psi}:(\overline{\widetilde{L}}^{\U(n)^{X}},t_p(\Delta)) \rightarrow  (\overline{L}^{\U(n)^{X}},t_p(\Delta)).$$

Using a compactness argument on the group $\U(n)$, it is not hard to verify that if $p,q\in \overline{\widetilde{L}}^{\U(n)^{X}}$ and
$\overline{\psi}(p)=\overline{\psi}(q)$ then $p$ and $q$ are equivalent.
Since $Orbit(p) = \{U^{-1}pU : U\in\U(n)\}$ has the cardinality of continuum $\frak c$ and
$\vert \overline{\widetilde{L}}^{\U(n)^X}\vert= \vert\beta L\vert= \vert\beta\omega\vert=2^{\mathfrak{c}},$

we obtain: $$2^{\mathfrak{c}}=\vert \overline{\widetilde{L}}^{\U(n)^X}\vert\leq \vert \overline{L}^{\U(n)^X}\vert\vert \U(n)\vert=\max\lbrace \vert \overline{L}^{\U(n)^X}\vert,\mathfrak{c}\rbrace.$$ Therefore $$\vert \overline{L}^{\U(n)^X}\vert\geq 2^{\mathfrak{c}}.$$

Applying \cite[Cor. 2.16]{Fer_Her_Tar_dichotomy}, it follows that $L$ contains a subset $P$ such that $\overline{P}^{\U(n)^X}$
is canonically homeomorphic to $\beta P$ (with $P$ equipped with the discrete topology). This completes the proof.
\end{enumerate}
\end{proof}

\bcor\label{resultado_7_ABEL}
Let $X$ be a \v{C}ech-complete abelian group.
If an infinite subset $G$ of $\widehat{X}$ is not equicontinuous, then $G$ contains a countably infinite $I_0$-set.
\ecor
\mkp

Next follows the proof of Theorem C.

\begin{proof}[\textbf{Proof of Theorem \ref{theorem_C}}]
If for every countable subset
$L\subseteq G$ and compact separable subset $Y\subseteq X$ we have that either
$\overline{L}^{K^Y}$ has countable tightness or $\vert\overline{L}^{K^Y}\vert \leq \mathfrak{c} $, then
$\overline{L}^{K^Y}$ may not contain any copy of $\beta\omega$.
By Theorem \ref{theorem_B}, this implies that $L\vert_Y$ is equicontinuous on $Y$.
Applying \cite[Theorem B]{Fer_Her_Tar2017}, it follows that $G$ is hereditarily equicontinuous on $X$, which implies that $G$ is equicontinuous
because $G$ consists of group homomorphisms.

\end{proof}
\bkp

\section{$I_0$-sets in abelian locally $k_{\omega}$ groups}

In this section, we study the existence of $I_0$-sets for abelian locally $k_{\omega}$ groups,
a large family of topological groups that includes, for example, all LCA groups,
the free abelian groups on a compact space and all countable direct sum of compact groups.
The proof of our main results are obtained using  methods of Pontryagin--van Kampen duality.
Therefore, we first recall some
basic defintions and facts about the Pontryagin duality of abelian groups.
From here on, all groups are supposed to be abelian and, therefore, we will use additive notation to deal with them.
In particular, we identify $\T$ with the additive group $[-1/2,1/2)$ by identifying $\pm 1/2$.

If $G$ is a topological group,
\emph{a character}\index{character!a character on a topological abelian group} on a topological abelian group $G$ is a continuous group homomorphism from $G$ to
the torus group $\T$.
The set of all characters on $G$, with pointwise addition, is a group.

For a topological abelian group $G$, let $\mathcal K(G)$ denote the family of all compact subsets of $G$.
For a set $A\subseteq G$ and a positive real $\epsilon$, define
$$[A,\epsilon]\defi \{\chi \in \widehat G : |\chi(a)|\leq\epsilon\mbox{ for all }a \in A\}.$$
The sets $[K,\epsilon]\subseteq\widehat G$, for $K\in\mathcal K(G)$ and $\epsilon>0$, form a neighborhood base at the trivial character,
defining the compact-open topology.
We write $\widehat G$ for the topological abelian group obtained in this manner.

A topological abelian group $G$ is \emph{reflexive} if the evaluation map $E$
$$E\colon G\to\widehat{\widehat G},$$
defined by $E(g)(\chi)=\chi(g)$ for all $g\in G$ and $\chi\in\widehat G$,
is a topological isomorphism.
By the Pontryagin--van Kampen theory, we know that
every LCA group is reflexive.
Furthermore, the dual of a compact group is discrete and the dual
of a discrete group is compact.
In general, the dual of a LCA group is also locally compact.
As a consequence, every compact abelian group is equipped with
the topology of pointwise convergence on its dual group.

\bdfn
Let $G$ be a topological abelian group.
For $A\subseteq G$, let $A^\rhd := [A,1/4]$\index{$A^\rhd$}.
Similarly, for $X\subseteq\widehat G$, let
$$X^\lhd :=\Bigl\{\, g\in G : |\chi(g)|\leq\frac{1}{4},\mbox{ for all }\chi \in X\,\Bigr\}.$$
\edfn
\mkp

The following facts are well known (see \cite{Banaszczyk1991}).

\blem\label{Upol}
For each neighborhood $U$ of $0$ in $G$, we have that $U^\rhd\in\mathcal K(\widehat G)$.
\elem

\bdfn
Let $G$ be a topological abelian group.
A set $A\subseteq G$ is \emph{quasiconvex} if $A^{\rhd\lhd}=A$.
The topological group
$G$ is \emph{locally quasiconvex} if it has a neighborhood base at the neutral element, consisting of quasiconvex sets.
$G$ is called \emph{Maximally Almost Periodic} (MAP, for short)
in the sense of von Neumann when $\widehat G$ separates the points in $G$,
which means that $G^+$ is Hausdorff and, as a consequence, the group $G$ is algebraically
injected in its Bohr compactification $bG$. Obviously, by definition,
every locally quasiconvex abelian group $G$ is MAP.
\edfn

For each set $A\subseteq G$, the set $A^\rhd$ is a quasiconvex subset of $\widehat G$.
Thus, the topological group $\widehat G$ is locally quasiconvex for each topological abelian group $G$.
Moreover, local quasiconvexity is hereditary for arbitrary subgroups.

The set $A^{\rhd\lhd}$ is the smallest closed, quasiconvex subset of $G$ containing $A$.

In the case where $G$ is a topological vector space, $G$ is locally quasiconvex in the present sense
if, and only if, $G$ is a locally convex topological vector space in the ordinary sense. 

If $G$ is locally quasiconvex, its characters separate points of $G$,
and thus the evaluation map $E\colon G\to \widehat{\widehat G}$ is injective.
For each quasiconvex neighborhood $U$ of $0$ in $G$, the set $U^\rhd$ is a compact subset of $\widehat G$
(Lemma \ref{Upol}), and thus $U^{\rhd\rhd}$ is a neighborhood of $0$ in $\widehat G$.
As $E[G]\cap U^{\rhd\rhd}=E[U^{\rhd\lhd}]=E[U]$, we have that $E$ is open \cite[Lemma 14.3]{Banaszczyk1991}.

The following theorem of Gl\"{o}ckner, Gramlich and Hartnick \cite{Glockner2010} states that
there exists a relation between the abelian locally $k_{\omega}$ groups and the abelian \v{C}ech-complete groups.

\bthm
If $G$ is an abelian locally $k_{\omega}$ group, the $\widehat{G}$ is \v{C}ech-complete.\\
Conversely, $\widehat{G}$ is locally $k_{\omega}$, for each abelian \v{C}ech-complete topological group $G$.
\ethm

Using this duality and Theorem \ref{theorem_B} we obtain:

\begin{proof}[\textbf{Proof of Theorem \ref{theorem_D}}]
Consider the abelian \v{C}ech-complete group $\widehat{G}$.
By means of the evaluation map
$E\colon G\to  \widehat{\widehat{G}}\subseteq C(\widehat{G},\T)$,
we can look at the sequence
$\lbrace g_n\rbrace_{n<\omega}$ as a subset of $C(\widehat{G},\T)$.
Furthermore, since  $\lbrace g_n\rbrace_{n<\omega}$ is not precompact in $G$,
it follows that $\lbrace g_n\rbrace_{n<\omega}$ is not equicontinuous on $\widehat G$.
Indeed, if it were equicontinuous on $\widehat{G}$, by Arzel\`{a}-Ascoli's theorem,
then $\lbrace g_n\rbrace_{n<\omega}$ would be precompact in $C_c(\widehat{G},\T)$,
the group  $C(\widehat{G},\T)$ equipped with the
compact open topology. Now, since $G$ is a locally quasiconvex $k$-space,
the evaluation map $E\colon G\to  \widehat{\widehat{G}}$
is a topological isomorphism onto its image (see \cite{Hernandez2001}).
Thus $\lbrace g_n\rbrace_{n<\omega}$ would also be precompact in $G$,
which is a contradiction.

Therefore, the sequence $\lbrace g_n\rbrace_{n<\omega}$ is not an equicontinuous set on $\widehat G$ and, by Theorem \ref{theorem_B},
contains an $I_0$-set.
\end{proof}
\mkp

The next result was proved in \cite[Lemma 4.11]{Galindo_Hernandez1999}.

\blem\label{Le_GH}
Let $G$ be a maximally almost periodic abelian group, $A$ a subset of $G$ and
let $N$ be a subset of $bG$  containing the neutral element such that
$A+N$ is compact in $bG$. If $F$ is an arbitrary subset
of $A$,  there exists $A_{0} \subseteq A$ with $|A_{0}|
\leq |N|$ such that
\[ cl_{bG}F \subseteq A_{0} +N + cl_{G^{+}}(F-F). \]
\elem
\mkp

We are now in position of proving Theorem \ref{theorem_E}.

\begin{proof}[\textbf{Proof of Theorem \ref{theorem_E}}]
Let $G$ be a locally quasiconvex, locally $k_\omega$ group and let $bG$ denote its Bohr compactification.

(i) Let $\mathcal P$ be a topological property implying functional boundedness and let $A$ be any subset of $G$ satisfying
$\mathcal P$ in $G^+$, which (by definition) is equipped with the weak topology generated by $\widehat G$;
that is to say $G^+\subseteq \T^{\widehat G}$.
Reasoning by contradiction assume that $A$ does not satisfy $\mathcal P$ in $G$. We claim that $A$
may not be a precompact subset of $G$. Indeed, if it were, since every locally $k_\omega$-group is complete
\cite[Remark 7.3]{Glockner2010}, it would follow that $\overline A^G$ would be compact in $G$. Therefore,
it would also be compact in $G^+$ that is equipped with a weaker topology.
Since any compact topology is (Hausdorff) minimal, this would imply that
the Bohr topology would coincide with the original topology of $G$ on the compact subset $\overline A^G$ and,
as a consequence, on $A$. Thus $A$ would have property $\mathcal P$ in $G$.

So, assume wlog that $A$ is not precompact in $G$.
As in the proof of Theorem \ref{theorem_D}, if we take
the abelian \v{C}ech-complete group $\widehat{G}$ and inject $G$ in $C(\widehat{G},\T)$
by means of the evaluation map
$E\colon G\to  \widehat{\widehat{G}}\subseteq C(\widehat{G},\T)$,
it follows that $A$ is not equicontinuous on $\widehat{G}$.
By \cite[Cor. 2.4]{Fer_Her_Tar2017}, it follows that there exists a countable subset
$F\subseteq A$ and a separable compact subset $X\subseteq \widehat G$
such that $F$ is not equicontinuous on $X$. Taking the closure in $\widehat G$ of the subgroup generated by $X$,
we may assume wlog that $X$ is a separable closed subgroup of $\widehat G$. On the other hand,
since $A$ is functionally bounded in $G^+$ and $X\subseteq \widehat G$, it follows that
$F$ is also functionally bounded in $G$, when the latter is equipped with weak topology generated by $X$.

Set $X^\bot\defi \{g\in G : \chi(g)=0\ \hbox{for all}\ \chi\in X\}$ and
take the quotient $G/X^\bot$, which clearly is a maximally almost periodic group
whose dual is $X$. Furthermore, the group $G/X^\bot$ is locally $k_\omega$ and
$\widehat{G/X^\bot}\cong X$, which is \v Cech-complete.
If $p\colon G\to G/X^\bot$ denotes the open quotient map and $bG/X^\bot$ denotes the Bohr compactification of
$G/X^\bot$, it follows that there is a canonical extension $p^b\colon bG\to bG/X^\bot$.
Therefore, we have that $p^b(F)$ is a functionally bounded subset of $p^b(G^+)=G^+/X^\bot$ that is not equicontinuous on $X$.
Indeed, if $p^b(F)$ were equicontinuous on $X$, then it would follow that $F$ would be equicontinuous on $X$,
which is not true. In other words, we may assume wlog that $A$ is a countable, functionally
bounded subset of $G^+$ that is not equicontinous on $\widehat G$, which is separable.

As in the proof of Theorem \ref{theorem_B}, by the \v{C}ech-completeness of $X$, there must be a
compact subgroup $C$ of $X$ such that $X/C$ is separable, complete and metrizable \cite{Brown1972}, thereby, a Polish space.

Let $A\vert_C\defi \lbrace g\vert_C:g\in A\rbrace\subseteq \widehat C$. We have two possible cases:
\begin{enumerate}[(1)]
\item $A_{\vert C}$ contains infinitely many different elements.
\item $A_{\vert C}$ only contains a finite number of different elements.
\end{enumerate}

\begin{enumerate}[-]
\item Case $(1)$: By Pontryagin duality, the dual of a compact group (equipped with the compact open topology) is discrete.
Therefore $A_{\vert C}$ is an infinite discrete subset of $\widehat C$ in the compact open topology on $C$,
which implies that $A_{\vert C}$ may not be equicontinuous on $C$. Applying Corollary \ref{resultado_5_COR_GT}, there is a nonempty subset $\Delta$
of $C$ and a countable subset $L$ of $A$ such that $L$ is separated by $\Delta$. Thus, by Lemma \ref{lem_i0}, $\overline L_{\vert \Delta}^{\T^{\Delta}}$
is canonically homeomorphic to $\beta L$ (where $L$ is equipped with the discrete topology), which yields
$$|\overline L_{\vert \Delta}^{\T^{\Delta}}|=2^\frak c.\qquad \qquad (\hbox{I})$$

On the other hand, by  Asanov and Velichko's  generalization of a well-known theorem of Grothendieck
about sets of continous functions defined on a compact space \cite{Asanov1981} (see also \cite[III.4.1.]{Arkhangelskij1992}),
we have that $C_p(\Delta,\T)$ (the space of continous functions of $\Delta$ into $\T$ equipped with the pointwise convergence topology)
is a $\mu$-space, which means that every closed functional bounded subset of $C_p(\Delta,\T)$ is compact.
Being $A$ a functionally bounded subset
of $G^+$, which is equipped with the pointwise convergence topology on $\widehat G\supseteq C\supseteq \Delta$,
it follows that $L_{\vert \Delta}\subseteq A_{\vert \Delta}$ is functionally bounded as a subset of  $C_p(\Delta,\T)$.
Thus $\overline L_{\vert \Delta}^{C_p(\Delta,\T)}$ must be compact and, as a consequence
$$\overline L_{\vert \Delta}^{C_p(\Delta,\T)}=\overline L_{\vert \Delta}^{\T^{\Delta}}.$$

However, by \cite[Lemma 3.4]{Galindo_Hernandez2004}, we have that $|\overline L_{\vert \Delta}^{C_p(\Delta,\T)}|\leq \frak c$, which is in contradiction with (I).

\item Case $(2)$: Set $H\defi \lbrace \phi_1,\ldots, \phi_m\rbrace\subseteq A$ such that for each $g\in A$ there is an element $\phi_{i}\in H$
satisfying that $g_{\vert C}=\phi_{i\vert C}$.
If we define $A_i=\lbrace g\in A:g_{\vert C}= \phi_{i\vert C}\rbrace$, then $A=A_1\bigcup\ldots\bigcup A_m$. Since $A$ is not equicontinuous,
there is $i\in\{1,\dots m\}$ such that $A_i$ is not equicontinuous. So, we may assume wlog that there is $g_0\in A$ such that
$g_{\vert C}= g_{0\vert C}$ for all $g\in A$. By Lemma \ref{resultado_6_1}, we know that $Ag_0^{-1}$ is not equicontinuous on $X$
and, by Lemma \ref{resultado_6_2}, it follows that $Ag_0^{-1}$ is  a $\frak B$-family on $X$.

Let $\pi_C:X\rightarrow X/C$ the canonical quotient map, which is open and continuous. Since $X/C$ is Polish and each ${g}g_0^{-1}$ factors through $X/C$,
we apply Theorem \ref{theorem_A} and Lemma \ref{lem_i0} in order to obtain $\Delta \subseteq X$ and ${L}\subseteq A$ such that
$$\overline{L_{\vert\Delta}}^{\T^{\Delta}}\simeq\overline{L_{\vert\Delta}g_0^{-1}}^{\T^{\Delta}}\simeq \beta {L},$$
which again yields
$$|\overline L_{\vert \Delta}^{\T^{\Delta}}|=2^\frak c.\qquad \qquad (\hbox{II})$$

On the other hand, repeating the same argument as in (1), we deduce that
$|\overline L_{\vert \Delta}^{\T^{\Delta}}|=|\overline L_{\vert \Delta}^{C_p(\Delta,\T)}|\leq \frak c$, which again is in contradiction with (II).
This completes the proof of (i).
\end{enumerate}

(ii) The proof of this part replays some of the steps followed to prove (i). For the reader's sake, we will avoid unnecessary
repetitions as much as possible.

Let $N$ be a closed metrizable subgroup of $bG$ and let $A$ be a subset $G$.
It is obvious that if $A+(N\cap G)$ is compact in $G$, then $b_N(A)$ is compact in $bG/N$.

In order to prove the non-trivial converse implication, again reasoning by contradiction,
assume that $b_N(A)$ is compact in $bG/N$ but $A+(N\cap G)$ is not compact in $G$. As $G$ is complete \cite[Remark 7.3]{Glockner2010},
this means that $A+(N\cap G)$, being closed in $G$, is not precompact in the topology inherited from $G$.
As in the proof of (i), it follows that there exists a countable subset
$F\subseteq A+(N\cap G)$ and a separable compact subset $X\subseteq \widehat G$
such that $F$ is not equicontinuous on $X$. Taking the closure in $\widehat G$ of the subgroup generated by $X$,
we may assume that $X$ is a separable closed subgroup of $\widehat G$.

Again, the quotient group $G/X^\bot$ is a MAP, locally $k_\omega$  group
whose dual group $X$ is \v Cech-complete.
If $p\colon G\to G/X^\bot$ denotes the open quotient map 
and $p^b\colon bG\to bG/X^\bot$ is the canonical extension to their Bohr compactifications,
we have that $p(A+(N\cap G))$ is contained in $p^b(A+N)=p(A)+p^b(N)$,
which is compact in $bG/X^\bot$. Applying Lemma \ref{Le_GH} to $p(F)$ and $p^b(N)$, we obtain
that there exists $A_{0} \subseteq p(A)$ with
$|A_{0}|\leq |p^b(N)|\leq \frak c$ such that
\[ cl_{bG/X^\bot}p(F) \subseteq A_{0} +p^b(N) + cl_{(G/X^\bot)^{+}}p(F-F). \]

Now, being the group $X$ separable, it follows that $G/X^\bot$ can be equipped with a metrizable
precompact topology. As a consequence $|G/X^\bot|\leq \frak c$. All in all, we obtain that
$|cl_{bG/X^\bot}p(F)|\leq\frak c$.

On the other hand, $p(F)$ is not equicontinuous as a subset of $C(X,\T)$ and, by Theorem \ref{theorem_D},
this means that it contains an $I_0$-set, which yields $|cl_{bG/X^\bot}p(F)|=|\beta\omega|=2^\frak c>\frak c.$
This is a contradiction that completes the proof.
\epf

\bkp

\bibliography{BiblioResultados}
\bibliographystyle{plain}

\bigskip

\end{document}